\documentclass[11pt]{amsart}
 \usepackage{graphicx}
\usepackage{xstring}
\usepackage{forloop}
\usepackage{bbm, amsmath, amsfonts, amscd, latexsym, amsthm, amssymb, graphicx, mathrsfs, comment}
\usepackage[colorlinks=true]{hyperref}

\usepackage{abstract}
\usepackage{mathtools}
\usepackage{caption}
\usepackage{MnSymbol}
\usepackage[dvipsnames]{xcolor}

\usepackage{tikz-cd}
\tikzset{
  symbol/.style={
    draw=none,
    every to/.append style={
      edge node={node [sloped, allow upside down, auto=false]{$#1$}}}
  }
}

\makeatletter
\newif\if@check@engine  \@check@enginetrue 
\makeatother
\usepackage[usenames,dvipsnames]{pstricks}

\usepackage{booktabs}

\newtheorem{theor}{\hspace{1cm}{\sc Theorem}}[section]
\newtheorem{utver}[theor]{\hspace{1cm}{\sc Proposition}}
\newtheorem{sledst}[theor]{\hspace{1cm}{\sc Corollary}}
\newtheorem{lemma}[theor]{\hspace{1cm}{\sc Lemma}}

\newtheorem*{utver*}{\hspace{1cm}{\sc Proposition}}
\theoremstyle{definition}

\newtheorem{defin}[theor]{\hspace{1cm}{\sc Definition}}
\newtheorem*{defin*}{\hspace{1cm}{\sc Definition}}
\newtheorem{exa}[theor]{\hspace{1cm}{\sc Example}}

\newtheorem{rem}[theor]{\hspace{1cm}{\sc Remark}}

\newtheorem{quest}[theor]{\hspace{1cm}{\sc Question}}

\newcommand{\ind}{\mathop{\rm ind}\nolimits}

\newcommand{\Vol}{\mathop{\rm Vol}\nolimits}

\newcommand{\codim}{\mathop{\rm codim}\nolimits}

\newcommand{\sing}{\mathop{\rm sing}\nolimits}

\newcommand{\conv}{\mathop{\rm conv}\nolimits}

\newcommand{\MV}{\mathop{\rm MV}\nolimits}

\newcommand*{\DashedArrow}[1][]{\mathbin{\tikz [baseline=-0.25ex,-latex, dashed,#1] \draw [#1] (0pt,0.5ex) -- (1.3em,0.5ex);}}

\newcounter{idx}

\newcommand{\rotraise}[1]{
  \StrLen{#1}[\slen]
  \forloop[-1]{idx}{\slen}{\value{idx}>0}{
    \StrChar{#1}{\value{idx}}[\crtLetter]
    \IfSubStr{tlQWERTZUIOPLKJHGFDSAYXCVBNM}{\crtLetter}
      {\raisebox{\depth}{\rotatebox{180}{\crtLetter}}}
      {\raisebox{1ex}{\rotatebox{180}{\crtLetter}}}}
}

\renewcommand{\emph}[1]{{\it {\color{NavyBlue} #1}}}

\def\R{\mathbb R}

\def\Z{\mathbb Z}

\def\C{\mathbb C}
\def\CC{({\mathbb C}^\star)}

\def\CP{\mathbb C\mathbb P}

\usepackage{xr-hyper}
\makeatletter

\newcommand*{\addFileDependency}[1]{%
\typeout{(#1)}%
\@addtofilelist{#1}
\IfFileExists{#1}{}{\typeout{No file #1.}}
}\makeatother

\newcommand*{\myexternaldocument}[1]{
\externaldocument[CR]{#1}
\addFileDependency{#1.tex}
\addFileDependency{#1.aux}
}

\myexternaldocument{crit}

\begin{document}

\begin{center}{\Large \sc Sch\"on complete intersections}

\vspace{3ex}

{\sc Alexander Esterov}
\end{center}

\begin{abstract}

A complete intersection $f_1=\cdots=f_k=0$ is sch\"on, if $f_1=\cdots=f_j=0$ defines a sch\"on subvariety of an algebraic torus for every $j\leqslant k$. 
This class includes nondegenerate complete intersections, critical loci of their coordinate projections, other simplest Thom–-Boardman and multiple point strata of such projections, generalized Calabi--Yau complete intersections, equaltions of polynomial optimization, hyperplane arrangement complements, and many other interesting special varieties.

We study their Euler characteristics, connectednes, Calabi--Yau-ness, tropicalizations, etc., extending (in part conjecturally) the respective classical results about nondegenerate complete intersections.
\end{abstract}

\tableofcontents

\section{Introduction}

Sch\"on subvarieties of an algebraic torus have a smooth toric compactification. Even before the notion was introduced by Tevelev \cite{tev}, such algebraic sets played an important role in geometry: a classical example is a nondegenerate complete intersection (abbreviated NCI) in the sense of Khovanskii (see \cite{khovcomp} or Deifintion \ref{defnondeg0} below).

One might have a false impression that nondegeneracy is close to the following.
\begin{defin*}
A complete intersection 
in the algebraic torus $\CC^n$ is said to be {\it sch\"on} (abbreviated to SCI) if, for every $i$, its first $i$ equations define 
a codimension $i$ sch\"on subvariety.
\end{defin*}
Our aim is to dispel this impression: many important complete intersections listed in the abstract are degenerate, yet sch\"on. %Among them are hyperplane arrangement complements and certain classes of algebraic sets naturally appearing e.g. in Galois theory (\cite{symm}),  enumerative and likelihood geometry (\cite{crit}).
NCIs are a key source of special varieties (such as Calabi--Yau and nearly rational ones, see e.g. \cite{ckp} for a demonstrative example). Passing to the much wider class of objects of similar nature may significantly enrich this source. 
This explains our motivation to extend classical results about NCIs to SCIs.

\subsection{Definition and examples}
\begin{defin}\label{deftropcomp}
1. A {\it complete intersection} in a torus $T\simeq\CC^n$ is a sequence of hypersurfaces $T=H_0\supset H_1\supset \cdots\supset H_k$ defined as the zero loci of regular functions $f_i:H_{i-1}\to\C$.

2. A toric compactification $X\supset T$ is called {\it tropical} for the complete intersection $\mathcal{H}=(H_i,f_i)_{i=1,\ldots,k}$, if the closure $\bar H_i\subset X$ intersects every orbit by its codimension $i$ subset.
\end{defin}

To describe $\bar H_2\subset X$ with equations $\varphi_1=\varphi_2=0$, we should take $\varphi_2$ to be a section of a line bundle on $\bar H_1=\{\varphi_1=0\}$ (rather than on $X$, as the default meaning of the term {\it toric complete intersection} would imply). We prefer to describe $H_i\subset T$ with polynomials and take closure.

\begin{defin}\label{defschon0}
The complete intersection $\mathcal{H}$ is said to be {\it sch\"on} (abbreviated as SCI), if the closure of every set $H_i$ at every its boundary point $x\notin T$ is smooth, reduced (as a component of the divisor $f_i=0$) and transversal to the orbit of $X$ passing through $x$.
\end{defin}
This property of $\mathcal{H}$ does not depend on the choice of a tropical compactification.

\begin{exa}\label{exascibkk0}
1. For generic Laurent polynomials $\Phi_i$ on $T\simeq\CC^n$ with prescribed Newton polytopes $P_i\subset M\simeq\Z^n$, the complete intersection $\Phi_1=\cdots=\Phi_k=0$ is sch\"on. 

2. {\it Critical complete intersections.} Let $\Phi$ be a generic polynomial with a prescribed Newton polytope, then the set $\Phi=0$ is smooth, and the critical locus of its coordinate projection $p$ is a Newtonian SCI. More generally, this is true for any system of equations involving $\Phi$ and its partial derivatives of arbitrary orders. 
We study such sets in \cite{crit}. 

3.{\it Qualitative polynomial optimization.} Given polynomials $f_0,\ldots,f_k$, the real extrema of $f_0$ on $f_1=\cdots=f_k=0$ are studied by polynomial optimization, and the number of complex extrema is known as the algebraic degree. The most studied cases include the maximum likelihood degree \cite{chks} ($f_0$ is a generic monomial), the polar degree \cite{dp} ($f_0$ is generic linear), and the euclidean distance degree \cite{dhost} ($f_0$ is the distance to a generic point). 
The algebraic degree is the number of critical points of the Lagrange multipliers function $F(x,\lambda)=f_0(x)+\lambda_1f_1(x)+\cdots+\lambda_kf_k(x)$. If $f_i$'s are generic polynomials with given supports, then $F$ is generic for its support, and results of \cite{crit} on critical complete intersections apply.

4. {\it All hyperplane arrangement complements} are SCI (as noticed in \cite{tev}). 

5. {\it Engineered complete intersections.} The general class of complete intersections covered in \cite{crit} includes (1-4) and many other interesting examples, such as  the following construction, allowing to engineer complete intersections with desired geometry.

6. A reflexive polytope decomposable into a Minkowski sum gives rise to a nondegenerate complete intersection Calabi--Yau, admitting a beautiful mirror construction \cite{bb94}. 
Let $r, g$ and $b$ be generic linear combinations of red, green and blue monomials in the unit cube, and let $f_1$ and $f_2$ be generic linear combinations of $r, g$ anf $b$. Then $f_1=f_2=0$ is a sch\"on CICY (an ellpitic curve) in $\CC^3$, corresponding to the \underline{non}-decomposable reflexive polytope $P$ on the right, in the sense explained below (Example \ref{exacy}). The basic geometry of this SCI (its connectedness, Euler characteristics, tropicalization) is covered by the results in \cite{crit}.

Note that $\{f_1=f_2=0\}$ is not a nondegenerate complete intersection, and not the zero locus of a pair of sections of line bundles on its tropical compactification, the $P$-toric variety.

\begin{center}

\includegraphics[scale=0.6]{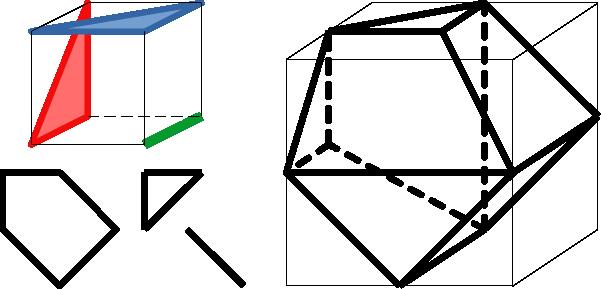}

{\it\small Figure 1: the reflexive polytope on the right is not a Minkowski sum of smaller lattice polytopes: 

the only decomposition of its projection along the green direction (on the bottom left) does not lift.}
\end{center}

\vspace{1ex}

7. {\it Generalized complete intersection Calabi--Yau varieties} (gCICY, \cite{gcicy}) belong to the class of sch\"on complete intersections, whose tropical compactification is a product of projective spaces. Though their reflexive polytope is a Minkowski sum of simplices, they exceed the class of nondegenerate CICYs that can be produced from this sum as in \cite{bb94}.

8. {\it Symmetric complete intersections.} Let $\Phi$ be a generic polynomial with a prescribed Newton polytope. Then the set $\Phi(x_1,x_2,x_3,\ldots,x_n)=\Phi(x_2,x_1,x_3,\ldots,x_n)=0$ 
splits into components that are Newtonian SCIs. We study such symmetric varieties in a separate paper \cite{symm}, motivated by their role in Galois theory and using the techniques developed here. 
\end{exa}

\subsection{Connectedness} Choose a smooth tropical compactification $X$ of a complete intersection $\mathcal{H}$ (it always exists). Every irreducibe component $\alpha$ of the {\it boundary} $\bar H_{i-1}\setminus T$ gives rise to:  

-- its ambient codimension 1 orbit of $X$, its local reduced defining equation $u=0$, and the corresponding linear surjection $l:M\to\Z$ of the character lattice $M\simeq\Z^n$ of the torus $T$ (the correspondence sends $l\in M^*$ to the orbit, whose cone in the fan of $X$ is generated by $l$).

-- 
the multiplicities $m_{i,\alpha}$ and $m$ of $\alpha$ as a component in the divisors of the functions $f_i$ and $u$, regarded as meromorphic functions on $\bar H_{i-1}$ (with positive multiplicity for poles);

\begin{defin}\label{defnewtpoyh} 1. The pairs $(l_{i,\alpha}:=m\cdot l,m_{i,\alpha})$ are called the {\it boundary data} of $\mathcal{H}$.

2. Its {\it Newton polyhedron} $N_i\subset M\otimes\R$ is defined by the inequalities $l_{i,\alpha}(x)\leqslant m_{i,\alpha}$. 
\end{defin}

Newton polyhedra non-trivially depend on the choice of defining equations $f_i$, not only the sets $H_i$, see Example \ref{exaempty}. 

We do not call them polytopes, because they may be unbounded.

Newton polyhedra do not depend on the choice of a tropical compactification. 
\begin{theor}\label{th0connect}

1. If $\mathcal{H}$ is SCI, then the set $H_k$ has at most isolated singularities.

2. If its Newton polyhedra 
have full dimension, then $\beta_j(H_k)=\beta_j(T)$ for the Betti numbers in all dimensions $j<n-k$. 
In particular, the set $H_k$ is connected unless it is a singular curve.

\end{theor}

The theorem is proved using Morse theory, see Section \ref{ssproofsconn}. (In the simplest case of $N_1=\cdots=N_k$ and smoe other additional assumptions it boils down to the classical Morse theoretic proof of the Lefschetz hyperplane section theorem for the respective toric variety.)

Note that the classical theorems about nondegenerate complete intersections (e.g. Theorem \ref{thkhconnect}) fit into this setting, when every hypersurface $H_i$ is defined in the previous one by a generic polynomial equation with a prescribed Newton polytope. 
Thus our Morse theoretic Theorem \ref{th0connect} gives a new proof to some of these classical results (see the end of Section \ref{ssproofsconn} for details). In the opposite direction, those classical results that are still not covered by Theorem \ref{th0connect}, suggest its possible extensions: see in particular Question \ref{qqqessent0} and the following one.
\begin{quest}\label{conjflat}
Do the components of a SCI in a torus $T$ always belong to pairwise distinct cosets of a certain subtorus?
\end{quest}
Results in \cite{symm} and \cite{kh15} (Theorem \ref{thkhconnect}), as well as many examples in the vein of Remark 1.5.3 
in \cite{crit}, substantiate this question (at least for the most interesting classes of SCIs).

\subsection{Newtonianity}
All the subsequent results in our paper address SCIs satisfying a mild extra assumption of being Newtonian. 
Theorem \ref{th0connect}.1 allows to verify this assumption (see Remark \ref{remnewtind}), and a positive answer to Question \ref{conjflat} would allow to drop it (see Remark \ref{remnonnewt}). This is one reason why we think this question matters. 
\begin{defin}\label{defnewtsci} A complete intersection $\mathcal{H}$ is said to be {\it Newtonian}, if in its boundary data $\R_+\cdot l_{i,\alpha}=\R_+\cdot l_{i,\beta}$ implies $m_{i,\alpha}\cdot l_{i,\beta}=m_{i,\beta}\cdot l_{i,\alpha}$ for all $i\leqslant k$. In this case, 
for any ray $l$ of the form $\R_+\cdot l_{i,\alpha}$, we correctly define the {\it Newton datum} $m_{i}^l:=m_{i,\alpha}/($lattice length of $l_{i,\alpha})$. 

\end{defin}
In other words, Newtonian complete intersections are those that can be nicely tropicalized: not only as varieties, but moreover as systems of equations, in the following sense. 

Recall that the {\it tropical fan} $F_k$ of the codimension $k$ set $H_k$ is the union of codimension $k$ cones of the fan $\Sigma$, whose corresponding orbits in the tropical compactification $X_\Sigma$ intersect the closure $\bar H_k$; the total multiplicity  of this intersection is called the {\it weight} of the cone in $F_k$.

\begin{defin}\label{tropeq}
The tropicalization of the equation $f_k$ of the Newtonian $\mathcal{H}$ is the continuous piecewise linear function $m_k:F_{k-1}\to\R$, defined on each codimension $k-1$ cone $S\in\Sigma$ as the linear function equal to $m_k^l$ on each of the $n-k+1$ generators $l\in S$.
\end{defin}
\begin{rem}
1. Tropicalizations do not depend on the choice of the tropical compactification $X_\Sigma$. This can be seen from their alternative definition (in case $\bar H_{k-1}$ is normal): a ray belongs to the fan $F_k\,\Leftrightarrow$ it is generated by the valuation $l$ of an analytic curve germ $\gamma:\C^*\to H_{k-1}$; in this case $m_k(l)$ equals the valuation of $f_k\circ\gamma$.

2. To make this definition work in the non-normal case, we should have let $\alpha$ 
in the definition of the boundary data $(l_{k,\alpha},m_{k,\alpha})$ to be a component of (the normalization of $\bar H_{k-1})\setminus H_{k-1}$, which would make the notion of a Newtonian complete intersection slightly more restrictive.

3. In particular, (1) implies that whether $\mathcal{H}$ is Newtonian does not depend on the choice of its smooth tropical compactification $X_\Sigma$. 
\end{rem}

\begin{exa}\label{exascibkk}
1. For generic Laurent polynomials $\Phi_i$ with prescribed Newton polytopes $P_i\subset M\simeq\Z^n$, the SCI $\Phi_1=\cdots=\Phi_k=0$ is Newtonian, its Newton datum $\tilde m^{\R_+\cdot l}_i$ for primitive $l\in M^*$ equals $\max l|_{P_i}$, so the tropicalization $m_i$ of $\Phi_i$ is the support function of $P_i$.

Its tropical compactification is a smooth toric variety $X_\Sigma$ whose fan $\Sigma$ is {\it compatible} with $P_i$'s (meaning that $m_i$'s are linear on every cone of $\Sigma$). 
Its Euler characteristics equals $$e_n(P_1,\ldots,P_k):=\frac{P_1}{P_1+1}\cdot\ldots\cdot \frac{P_k}{P_k+1}.\eqno{\rm (BKK)}$$Recall how to evaluate this classical BKK formula from \cite{kh75}: one should expand the rational function $P_i/(P_i+1)$ as $P_i-P_i^2+P_i^3-\cdots$, open the brackets, and evaluate every product of $n$ polytopes as their mixed volume (Definition \ref{defpolyh}), omitting the other products.
Note that the Newton polyhedron $N_i$ for this SCI may strictly contain $P_i$, but plugging it instead of $P_i$ in (BKK) yields the same number (the Euler characteristics). 

Below we extend this formula to all Newtonian SCIs.
To the best of our knowledge, the Betti numbers of nondegenerate complete intersections are not completely computed in the literature; all the more so we do not achieve this for arbitrary Newtonian SCIs.

2. All other SCIs in Example \ref{exascibkk0} are Newtonian as well. In particular, since Betti numbers of a hyperplane arrangement complement are known \cite{os}, it may serve as a valuable sandbox in the study of Betti numbers of arbitrary Newtonian SCIs.

\end{exa}
\begin{rem}

Proving our main results, as a byproduct we obtain two proofs of formula (BKK) different from the original paper \cite{kh75}: a Morse theoretic one (sketched in Section \ref{ssproofsconn}), and a deformation theoretic one (sketched in Section \ref{sproof2}). 
\end{rem}

\subsection{Tropicalization}

\begin{defin}\label{deftci}
A {\it tropical complete intersection} (TCI) is a sequence of tropical fans $F_k,\,\codim F_k=k$, in $F_0:=\R^n$, and continuous piecewise linear functions $m_{k}:F_{k-1}\to\R$ such that $F_{k}$ is the corner locus \cite{mikhicm}, or the Weil divisor \cite{rau07} of $m_k$ (denoted by $\delta(m_k\cdot F_{k-1})$).

\end{defin}

\begin{rem}\label{exaropnewt0}
For a Newtonian $H_k\subset\cdots\subset H_1\subset T\simeq\CC^n$, the tropical fans $F_i$ of the sets $H_i$ and the tropicalizations $m_i$ of their defining equations form a TCI, so that
$$F_k=\delta(m_k\cdot\delta(m_{k-1}\cdot\ldots\cdot\delta(m_2\cdot \delta m_1)\ldots)).\eqno{(*)}$$\end{rem}\noindent {\it{}Proof} proceeds by induction on $k$: for any variety $V\subset T$ and a regular function $f:V\to\C$ with tropicalizations $F$ and $m:F\to\R$, the tropicalization of the divisor $f=0$ is $\delta(m\cdot F)$.\hfill$\square$
\begin{exa}\label{exaropnewt}
In particular for $k=n$, the 0-dimensional fan $F_n$ consists of the cone $\{0\}$ with multiplicity $|f_n=0|$, i.e. formula $(*)$ counts the points (with multiplicities) of the 0-dimensional set $H_n=\{f_n=0\}$ in terms of the tropicalizations $m_i$. Especially for a nondegenerate complete intersection $\Phi=0$ of Example \ref{exascibkk}.1, the tropicalizations $m_i$ are the support functions of the Newton polytopes $P_i$, so the right hand side of $(*)$ becomes the mixed volume of the Newton polytopes (as shown combinatorially in \cite{mmjpoly}), and thus the equality $(*)$ becomes the Bernstein-Kouchnirenko formula \cite{bernst}: $|\Phi=0|=\MV(P_\bullet)$. 

How equality $(*)$ works for degenerate complete intersections is shown in Example \ref{exaempty}. 
\end{exa}
\begin{exa}\label{exacy} For a Newtonian SCI $H_k\subset\cdots\subset H_1\subset T$, assume that the function $\sum_i m_i$ extends from 
the fan $F_{k-1}$ to the support function $m:M^*\to\Z$ of a convex polytope $P$, and $F_{k}$ is a union of certain cones $S_i$ from the dual fan of $P$, each of whom as a semigroup $S_i\cap M^*$ is generated by lattice points on the boundary of the dual polytope $\{a\,|\, m(a)=1\}$.

Then, by the adjunction formula, the respective toric compactification of $H_k$ is Calabi--Yau.
\end{exa}
\begin{quest}
This in particular happens if $P$ is reflexive, as in Example \ref{exascibkk0}.6 above.
If, moreover, all functions $m_i$ are convex (i.e. are themselves support functions of lattice polytopes, which then sum up to $P$), then we have the Batyrev-Borisov construction of mirror CYCI \cite{bb94}. Dropping the convexity assumption, but assuming that $P$ is a Minkowski sum of simplices, so that $H_k$ is gCYCI, there are great recent advances in the same direction (see \cite{bh} for the big picture).
It would be interesting to know to what extent it generalizes to Newtonian sch\"on, and, especially, to engineered complete intersections (in the sense of \cite{crit}).
\end{quest}
\subsection{Topology}
Theorem \ref{th0connect}.2 is not a combinatorial expression for the topological Betti numbers of a SCI. Such an expression is given below for Newtonian SCIs.
One might expect that the Euler characteristics is expressed in terms of its Newton polyhedra (Definition \ref{defnewtpoyh}) by the same classical formula as for non-degenerate complete intersections (Example \ref{exascibkk}.1).

This indeed holds true if an SCI can be defined by equations, whose tropicalizations 
are convex (and thus equal the support functions of the Newton polyhedra). However, this fails for the most interesting SCIs (such as the ones in Example \ref{exascibkk0}).
\begin{theor}[proved in Section \ref{sec5}]\label{th0tropeuler}
The diffeomorphic type of a Newtonian SCI $H_k\subset\cdots\subset H_1$ does not change under deformations preserving its tropicaliztion $(F_i,m_i)$ and singularities (recall that they are ICIS, if any). The Euler characteristics of $H_k$ equals
$$\frac{\delta m_k}{1+\delta m_k}\cdot\ldots\cdot\frac{\delta m_1}{1+\delta m_1}+(-1)^{n-k}(\mbox{the sum of the Milnor numbers of singularities of }H_k).\eqno{(**)}$$

\end{theor}

\begin{rem}\label{remtropeuler}
1. Formula $(**)$ 
should be evaluated as follows (cf Example \ref{exascibkk}.1): 

-- expand every multiplier $\delta m_k/(1+\delta m_k)$ as $\delta m_i-(\delta m_i)^2+(\delta m_i)^3-\ldots$;

-- open the brackets, and evaluate every resulting monomial as an iterated corner locus $(*)$;

-- if the resulting tropical fan is not 0-dimensional, then omit it;

-- the rest is a sum of 0-dimensional tropical fans, i.e. the cone $\{0\}\subset\R^n$ with a multiplicity. This multiplicty (plus the Milnor numbers, if any, see e.g. \cite{lij}) is the Euler characteristics.

2. If we do not omit positive dimensional fans when evaluating formula $(**)$, then what we get is the full tropical characteristic class of $H_k$ (as defined in \cite{e13}): its 0-dimensional component is the Euler characteristics (recovering Theorem \ref{th0tropeuler}), and its highest dimension component is the tropicalization (recovering $(*)$).

In particular, descending (as in \cite{fs}) from the ring of tropical fans to cohomology ring of a toic variety $X_\Sigma\supset T$ for a sufficiently fine fan $\Sigma$, the tropical characteristic class $\langle S\rangle$ evaluates to the CSM class of the closure of $S$ in $X_\Sigma$. (A sufficient condition on $\Sigma$ is that the tropical fan  representing each component of the characteristic class $\langle S\rangle$ is a union of cones of $\Sigma$.)

3. Using the notion of tropical fans with polynomial weights \cite{mmjpoly}, we can rewrite $(**)$ as $$\frac{\delta^{n-k}}{(n-k)!}\left(\frac{1}{(1+m_1)\cdots(1+m_k)}\cdot\delta m_k\cdot \ldots\cdot\delta m_1\right);\mbox{ or, further, } \frac{\delta^n}{n!} \left(\frac{m_1}{1+m_1}\cdots\frac{m_k}{1+m_k}\right),$$ arbitrarily extending $m_i$ to continuous piecewise linear functions on $\R^n$. In practice, it is more computationally effective. Another way to compute such expressions is given in \cite{kppoly}.

4. If a Newtonian SCI $H_k$ has full dimensional Newton polyhedra, then we can find its Betti numbers as follows: $\beta_i(H_k)=0$ for $i>\dim H_k$ because $H_k$ is an affine variety, $\beta_i(H_k)=\beta_i(T)$ for $i<\dim H_k$ by Theorem \ref{th0connect}, and $\sum_i (-1)^i \beta_i(H_k)$ is computed in Theorem \ref{th0tropeuler}.

5. It would be interesting to further investigate geometry of SCIs in terms of their tropicalizations, especially the Hodge-Deligne numbers. %because of the possibility to construct new Fano and Calabi-Yau varieties as Newtonian SCIs (see Example \ref{exacy}). 
Recent advances for gCYCIs suggest methods and support the expectation that the result is indeed determined by the tropicalization. 

6. Note that topology of a Newtonian SCI is not determined by its tropicaliztion, because, even for line arrangement complements in the plane, it is not defined by the respective matroid \cite{rybn}. However, the difference is seen only at the level of the fundamental group, while Betti \cite{os} and Hodge-Deligne numbers \cite{bs} of hyperplane arrangement complements are determined by their tropicalization. The same is expectable for arbitrary Newtonian SCIs.
\end{rem}
\subsection{The structure}

In Section \ref{sec6} we prove that a deformation of a sch\"on variety does not change its topology, as long as its tropical fan does not change (Proposition \ref{schontopol}). Then we show how this applies to nondegenerate complete intersections. Along the way, we recall the underlying notions such as the tropical fan and the nondegenerate complete intersection.

Section \ref{smorse} introduces Morse theoretic techniques proving Theorem \ref{th0connect}. 
Section \ref{sec5} presents Theorem \ref{th0euler}, recursively computing the Euler characteristics of SCIs. It is more flexible than its corollary, Theorem \ref{th0tropeuler}. 
Finally, in Section \ref{sdiscus} we highlight some open questions about SCIs and tropical complete intersections that we consider important.

\section{Tropicalizations and topology}\label{sec6}

\subsection{Topology of sch\"{o}n varieties} We prove that the topology of a sch\"on variety is preserved under deformations as long as its tropical fan does not change. See e.g. \cite{ms} as a general referens on tropical geometry.
\begin{defin}
A smooth toric variety $X_\Sigma$ compactifying $T\simeq\CC^n$ is said to be a tropical compactification
of an algebraic set $V\subset T$ of pure codimension $k$, if the closure $\bar V$ in $X_\Sigma$ does not intersect orbits of dimension smaller than $k$.
\end{defin}
\begin{rem}
1. Every algebraic set has a tropical compactification.

2. The definition implies that $\bar V$ intersects every orbit of dimension $d$ by a set of dimension $d-k$ (i.e. properly). In particular, it intersects $k$-dimensional orbits by finite sets, and the following definition makes sense.
\end{rem}
\begin{defin}
If $X_\Sigma$ is a tropical compactification of $V$, then the tropical fan of $V$ is defined as the union of codimension $k$ cones $C$ in the fan $\Sigma$ endowed with the weights defined as the intersection indices
$$\bar V\circ(\mbox{the orbit corresponding to } C).$$
\end{defin}
The definition does not depend on the choice of a tropical compactification.

\begin{utver}\label{schontopol}
Let $V_u\subset T$ 
be a family of varieties  parameterized by the points $u$ of a connected variety $U$ and having the same topical fan $S$. 
If the closure of every $V_u$ in a toric variety $X_\Sigma\supset T$ is smooth and transversal to the orbits, then all of $V_u$'s are diffeomorphic.
\end{utver}

\begin{rem} 1. Slightly more generally, assume that each $\bar V_u$ is smooth near $X_\Sigma\setminus T$ and transversal to every orbit there. Then $V_u$ has only isolated singularities.

2. In this case, if we additionally assume that the number and type of singularities does not depend on $u$, then neither does the topological type of $V_u$.
\end{rem}
\begin{proof}
Choose a smooth fan $\Sigma$ including all the cones of $S$, then $X_\Sigma$ is a tropical compactification of every member of the family. With no loss in generality, we can assume that $U$ is a neighborhood of $0\in\C$. 

The closure $\bar V$ intersects the 0-fiber of the projection $X_\Sigma\times U\to U$ by a codimension $k$ set $W$ which contains $\bar V_0$ and a priory may contain some other component $V'\in(X_\Sigma\setminus T)\times\{0\}$.

We shall prove that such $V'$ does not exist, using the fact that, for any $k$-dimensional orbit $\mathcal {O}\subset X_\Sigma$, the intersection index of $\bar V$ with $\mathcal {O}\times\{\varepsilon\}$ does not depend on $\varepsilon$; denote it by $N_{\mathcal{O}}$. Since for all but finitely many $\varepsilon$ the set $\bar V$ intersects $X_\Sigma\times\{\varepsilon\}$ exactly by $\bar V_\varepsilon$ (by considerations of dimension), the aforementioned intersection index $N_{\mathcal{O}}$ is the weight of the cone corresponding to the orbit $\mathcal{O}$ in the tropical fan $S$.

If, toward the contradiction, $W$ is the union of $V_0$ and $V'\ne\varnothing$, then the intersection number of $V'$ with some $k$-dimensional orbit would be positive, thus the intersection number of this orbit $\mathcal{O}$ with $V_0$ would be strictly smaller than $N_{\mathcal{O}}$, which would distinguish the tropical fan of $V_0$ from $S$.

Now we know that $\bar V\cap(X_\Sigma\times\{0\})=\bar V_0$, and thus a sufficiently close fiber $\bar V_\varepsilon$ is contained in a small metric neighborhood of $\bar V_0$. Introduce in a small neighborhood of $\bar V_0$ a metric such that $\bar V$ is perpendicular to $\mathcal{O}\times U$ for every orbit $\mathcal {O}\subset X_\Sigma$: near any given point it can be the standard metric in a local system of coordinates in which the normal crossing divisor $D:=(X_\Sigma\setminus T)\times U$ and $\bar V$ consist of coordinate planes; then such local metrics can be glued into a global one with a partition of unity.

Shifting $\bar V_0$ perpendicularly to itself in the sense of this metrics gives its diffeomorphism with $\bar V_\varepsilon$ preserving its intersection with the divisor $D$, and thus inducing the diffeomorphism between $V_0$ and $V_\varepsilon$.
\end{proof}

\subsection{Nondegenerate complete intersections} We recall the notion of nondegenerate complete intersection, and how the material of the preceding subsection applies to it. 

Given a finite set $A$ in the character lattice $M\simeq\Z^n$ of the torus $T\simeq\CC^n$, the vector space generated by characters $m\in A$ as functions on $T$ is denoted by $\C^A$. 
An element $f\in\C^A$ is regarded as a function on $T$, in coordinates $T\simeq\CC^n$ it is given by a Laurent polynomial.

\begin{defin}\label{defrestr}
1. Given a finite set $B$ in the character lattice $M$, and a Laurent polynomial $f(x)=\sum_{m\in M} c_m m(x)$, the restriction of $f$ to $B$ is defined as $f|_B(x):=\sum_{m\in B} c_m m(x)$. 

2. For a linear function $\gamma:M\to\Z$, the {\it $\gamma$-support face} $A^\gamma$ of $A\subset M$ is the set of all points at which $\gamma$ attains its maximum on $A$. The {\it $\gamma$-leading component} of $f\in\C^{A}$ is defined as $f|_{A^\gamma}(x)$.
\end{defin}
\begin{rem}
Note that the notion of the $\gamma$-leading component is relative to the choice of the support set $A$: for instance, if $n=1$, the 1-leading component of $px+q$ as an element of $\C^{\{0,1\}}$ equals $px$, but becomes 0 once we regard $px+q$ as an element of $\C^{\{0,1,2\}}$.
\end{rem}

\begin{defin}\label{defnondeg0}
A collection of polynomials $f=(f_1,\ldots,f_k)$
supported at finite sets $A_1,\ldots,A_k\subset M$ is said to be {\it nondegenerate} with respect to these sets, if, for every linear function $\gamma:M\to\Z$, the system of equations $f_1|_{A_1^\gamma}=\cdots=f_k|_{A_k^\gamma}=0$ is regular (i.e. the differentials of the equations are linearly independent at every root of the system). 
\end{defin}
For any given $\gamma$, the set of $f\in\C^{\mathcal A}=\C^{A_1}\oplus\cdots\oplus\C^{A_k}$ satisfying this condition is contained in a proper algebraic subset $\Sigma_\gamma\subset\C^{\mathcal A}$. Importantly, this set is the same for all linear functions $\gamma$ with the same support faces $A_i^\gamma$. 
Since such collections of support faces are finite in number, we have finitely may sets of the form $\Sigma_\gamma$ for the infinitely many linear functions $\gamma$. Thus the union $\bigcup_\gamma\Sigma_\gamma$ of finitely many proper algebraic subsets is itself a proper algebraic subset (and moreover a hypersurface \cite{adv}), i.e. we have the following. 
\begin{theor}[\cite{khovcomp}] \label{thbkknondeg}
Nondegenerate collections form a non-empty Zariski open subset in the space $\C^{\mathcal A}:=\C^{A_1}\oplus\cdots\oplus\C^{A_k}$.

\end{theor}
\begin{rem}\label{remnondeg}

1. Tuples $f$ outside the set $\bigcup_\gamma\Sigma_\gamma$ are nondegenerate. In particular, their zero loci $\{f=0\}$ are  smooth: 
a system $f=0$ is by definition regular unless $f\in \Sigma_0$.

2. Tuples $f$ outside of the smaller set $\bigcup_{\gamma\ne 0}\Sigma_\gamma$ are called {\it nondegenerate at infinity}. The sets $f=0$ defined by them are not necessarily smooth, but have at most isolated 
singularities. 

3. 
A nondegenerate system of equations $f_1=\cdots=f_k=0$ may be not a SCI: to ensure that it is a SCI we should assume than $f_1=\cdots=f_i=0$ is nondegenerate for every $i\leqslant k$. 
\end{rem}

\subsection{Tropicalization and topology of nondegenerate complete intersections}
Geometric characteristics of a nondegenerate complete intersection can be expressed in terms of support sets of its equations. 
In particular, the Euler characteristics, the number of components and the tropical fan of non-degenerate $f=0,\,f\in\C^{\mathcal{A}}$, are expressed in terms of mixed volumes of $A_i$'s and their faces.
For the Euler characteristics, see Example \ref{exaropnewt}. Below we recall the expressions for the other two. 

\begin{defin}\label{defpolyh} The {\it lattice mixed volume} $\MV$ is the (unique) real-valued function of $n$ convex sets in $\R^n$, which is 
i) symmetric in its arguments, 

ii) linear in each argument w.r.t the Minkowski summation $B_1+B_2:=\{b_1+b_2\,|\,b_i\in B_i\}$,

iii) has the diagonal value $\MV(B,\ldots,B)$ equal to the lattice volume $\Vol_\Z (B):=n!\Vol (B)$. \end{defin}

\begin{rem}\label{remmv}
1. The multiplier $n!$ assures that the mixed volume of lattice polytopes (i.e. polytopes whose vertices belong to $\Z^n$) is integer. When we refer to mixed volumes of subsets of $\Z^n$, we imply mixed volumes of their convex hulls in $\R^n$.

2. More generally, for any lattice $M\simeq\Z^n$ and $n$ convex sets in its ambient space $\R\otimes M$, their $M$-mixed volume is defined as the mixed volume of their images under an induced isomorphism $\R\otimes M\xrightarrow{\sim}\R\otimes\Z^n=\R^n$ (the value does not depend on the choice of the identification).

3. More generally yet, if $k$ convex sets can be shifted to the same $k$-dimensional rational subspace $K\subset \R\otimes M$, their lattice mixed volume can be defined as the $K\cap M$-mixed volume of their shifted copies in $K$. The value does not depend on the choice of the shifts (and is 0 if $K$ is not unique, i.e. if the sets can be shifted to a less than $k$ dimensional plane). 
\end{rem}
We refer to \cite{ew} for a detailed introduction to mixed volumes in the presence of lattice.

\begin{theor}[\cite{kh15}]\label{thkhconnect}
The number of components of non-degenerate $f_1=\cdots=f_m=0,f\in\C^{\mathcal{A}}$, is the mixed volume of the maximal proper subtuple of $\mathcal{A}$, for which it is defined in the sense of Remark \ref{remmv}.3. (The mixed volume of the empty tuple is defined to be 1 here.)
\end{theor}
\begin{rem}
1. One can easily see that the number of components is not smaller than that. Assuming wlog that the maximal such subtuple is $A_1,\ldots,A_k$, we can choose coordinates in the character lattice $M\simeq\Z^n$ and in the torus $T\simeq\CC^n$ so that $A_1,\ldots,A_k$ can be shifted to the first $k$-dimensional coordinate plane. Thus $f_1,\ldots,f_k$, up to monomial multipliers, depend only on the first $k$ variables, and thus have $\MV(A_1,\ldots,A_k)$ roots in the torus $\CC^k$. 

The initial set $\{f=0\}$ is in the preimage of these roots under the projection $\CC^n\to\CC^k$. The theorem claims that the preimage of each root contains exactly one component of $\{f=0\}$.

2. As a corollary, a nondegenerate complete intersection is connected, if all of its Newton polytopes are full dimensional: this is a special case of Theorem \ref{th0connect} as well. More generally, a nondegenerate complete intersection is connected, unless its tropical fan is contained in a hyperplane. We now recall the description of the tropical fan.
\end{rem}

\begin{utver}\label{ptropnondeg} Assume that $f\in\C^{\mathcal A}$ is nondegenerate at infinity.

1. The set $\{f=0\}$ has at most isolated singularities. Its tropical fan $F$ consists of all $l\in M^*$ such that $A_1^l,\ldots,A_k^l$ cannot be shifted to the same $(k-1)$-dimensional sublattice of $M$.

2. For generic $l\in F$, they can be shifted to a $k$-dimensional sublattice $L$, and the multiplicity of the fan $F$ at this point $l$ equals the lattice mixed volume of $A_1^l,\ldots,A_k^l$ in $L$.
\end{utver}

\begin{sledst}\label{bkktopol}
All nondegenerate $f=0,\, f\in\C^{\mathcal A}$, 
have the same topological type.
\end{sledst}
\begin{proof}
By the preceding proposition, the sets $f=0$ for all nondegenerate $f$ have the same tropical fan. Thus Proposition \ref{schontopol} applies to the family of all such complete intersections.
\end{proof}
\begin{rem}
We put it here just to illustrate Proposition \ref{schontopol}. A stronger fact is proved in \cite{e13}: all nondegenerate complete intersections in $\C^{\mathcal A}$ have the same topological type, and nondegeneracy is the weakest such condition of general position. More precisely, almost all complete intersections that are not non-degenerate have a different topological type. 
\end{rem}

\section{Morse theory for Laurent polynomials}\label{smorse}
\subsection{Stratified Morse theory} Given a smooth function on a compact manifold $\varphi:M\to\R$, denote its lower level set $\{x\in M\,|\,\varphi(x)\leqslant c\}$ by $\{\varphi\leqslant c\}$ or $M_{\leqslant c}$. 
\begin{theor}[Morse toolkit] 1. If $c$ is the value of $\varphi$ at finitely many critical points $x_\alpha$, and all of them are Morse, then $\{\varphi\leqslant c\}$ is homotopy equivalent to $\{\varphi<c\}$ with a cell of dimension $\ind x_\alpha$ attached for every $\alpha$. In particular, if $c$ is regular, then $\{\varphi\leqslant c\}\sim\{\varphi<c\}$.

2. If values in $[c,c')$ are regular, then $\{\varphi< c\}$ is a deformation retract of $\{\varphi< c'\}$.

3. If values in $(c,c']$ are regular, then $\{\varphi\leqslant c\}$ is a deformation retract of $\{\varphi\leqslant c'\}$.
\end{theor}
This allows to trace the changes in topology of $\{\varphi\leqslant c\}$ as $c$ increases. We want to extend it to the following setting.

\begin{defin}\label{morsish}
Let $V$ be a smooth compact complex algebraic variety with a normal crossing divisor $D\subset V$. A meromorphic function $f:V\DashedArrow\CP^1$ is 
{\it Morse on} $(V,D)$, if:

1) $D$ and $\{f=0\}$ have no common components, and $D\supset\{f=\infty\}$;

2) 
critical points of $f$ outside of $D\cup\overline{\{f=0\}}$ are isolated;

3) the restriction of $f$ to every stratum of $D$ has no critical values besides $0$ and $\infty$; 

4) the closure $\overline{\{f=0\}}$ is transversal to every stratum of $D$ at their every common point.

\end{defin}

The conditions of Definition \ref{morsish} are designed to ensure that we have the full Morse toolkit for $\varphi=|f|$ on $U$, despite the latter is not compact.
\begin{rem}\label{realmilnor}
1. Note that the critical points of $\varphi$ and $f$ outside of $D\cup\overline{\{f=0\}}$ are the same (thus isolated). 
By the Milnor number of a critical point $x_\alpha$ of $\varphi$ we mean that of $f$. 

This number $\mu_\alpha$ can be indeed read off from $\varphi$ rather than $f$: for a small neighborhood $U\ni x_\alpha$, 
the set $\{\varphi\leqslant\varphi(x_\alpha)\}\cap U$ is homotopy equivalent to $\{\varphi<\varphi(x_\alpha)\}\cap U$ with $\mu_\alpha$ cells of real dimension $\dim_\C U$ attached.

2. Importantly, condition 2 in the definition is redundant, if the compliment to $\{f=0\}\cup D$ is affine. This will be proved together with the following.
\end{rem}

\begin{theor}\label{toricmorse} Let $f$ be a Morse meromorphic function on $(V,D)$, set $U:=V\setminus D$.

1. If $c\ne 0$ is the value of $\varphi:=|f|$ at its critical points $x_\alpha\in U$ with Milnor numbers $\mu_\alpha$, then the set $\{\varphi\leqslant c\}\cap U$ is homotopy equivalent to $\{\varphi<c\}\cap U$ with $\sum_\alpha\mu_\alpha$ cells of real dimension $\dim_\C U$ attached. In particular, if $c$ is a regular value, then $\{\varphi\leqslant c\}\cap U\sim\{\varphi<c\}\cap U$.

2. If values in $[c,c')\subset(0,+\infty)$ (with possibly $c'=+\infty$) are regular for $\varphi$ on $U$, then the set $\{\varphi< c\}\cap U$ is a deformation retract of $\{\varphi< c'\}\cap U$.

3. If values in $(c,c']\subset(0,+\infty)$ (with possibly $c=0$) are regular for $\varphi$ on $U$, then the set $\{\varphi\leqslant c\}\cap U$ is a deformation retract of $\{\varphi\leqslant c'\}\cap U$.
\end{theor}
\begin{exa}
1. Note that statements (2) and (3) obviously fail without condition (3) in Definition \ref{morsish}: look at $f(x:y:z)=y/z$ on $V=\CP^2$ with $D=\{xz=0\}$ as a positive example, and the same with $f(x:y:z)=(y/z+x/z+1)(y/z-x/z+1)+1$ as a negative one. 

2. Note that statement (3) for $c=0$ (unlike (4) for $c'=\infty$) obviously fails without condition (3) in Definition \ref{morsish}: look at $f(x:y:z)=y/z$ on $V=\CP^2$ with $D=\{xz=0\}$ as a positive example, and the same with $D=\{x(x+y)z=0\}$ as a negative one. 
\end{exa}
\begin{rem}\label{remmorseicis}
If we allow $V$ to have isolated complete intersection singularities outside $D$, with the convention that each of them is regarded as a critical point of $\varphi$, then the Milnor number at such point is defined as in Remark \ref{realmilnor}, and Theorem \ref{toricmorse} and its proof hold true.
\end{rem}
We prove Theorem \ref{toricmorse} at the end of this section, using the following auxiliary construction.

\subsection{Resolving the indeterminacy locus} Given a meromorphic function on a variety $V$ (i.e. a regular map $\varphi:U\to\CP^1$ on a Zariski open subset $U\ne\varnothing$ of $V$), we can extend $\varphi$ by regularity to the maximal possible $U$, and the complement $V\setminus U$ is called the {\it indeterminacy locus} of $\varphi$. It can be resolved in the following sense.
\begin{utver} There exists a degree 1 proper map $p:\tilde V\to V$ such that the indeterminacy locus of $\tilde\varphi:=\varphi\circ p$ is empty.
\end{utver}
\begin{proof} $\tilde V$ is the resolution of singularities for the closure of the graph of $\varphi$ in $V\times\CP^1$.\end{proof}
However, this one is not convenient for our purpose, because in general $p:\{\tilde\varphi=0\}\to V$ is not a homeomorphism onto $\{\tilde\varphi=0\}$, by the following reason.
\begin{exa}\label{examerom}
Let $\varphi(x)=x_1^{d_1}/(x_2^{d_2}\cdots x_m^{d_m})$ on $\C^m$ with $d_i>0$. Define the fan $\Sigma$ in $\Z^m$ as the common subdivision of the following two:

-- the positive orthant $C$ together with all of its faces;

-- the plane $d_1u_1=d_2u_2+\cdots+d_mu_m$, together with its two adjacent half-spaces $H_+$ and $H_-$.

The fan $\Sigma$ has a smooth subdivision $\tilde\Sigma$, and the corresponding projection of toric varieties $\tilde p:X_{\tilde\Sigma}\to X_C\simeq\C^m$ is the sought resolution: the pullback $\tilde \Phi:=\varphi\circ\tilde  p$ has no indeterminacy points. However, note that for most of $(d_1,\ldots,d_m)$, every $\tilde\Sigma$ has one-dimensional cones $R_\alpha$ in the interior of $H_+$ and $H_-$, different from the edges of $C$. 
The $R_\alpha$-orbits of $X_{\tilde\Sigma}$ will be the components of the divisors $\{\tilde \Phi=0\}$ and $\{\tilde\Phi=\infty\}$ respectively, which are mapped by the projection $\tilde p$ to the indeterminacy locus of $\varphi$ (not injectively).

We instead wish to consider a singular resolution of the indeterminacy locus, free from this drawback, such as $p:X_\Sigma\to X_C\simeq\C^m$ in the preceding example.

The function $\Phi:=\varphi\circ p$ is not Morse in the sense of Definition \ref{morsish}, because its domain is not smooth. However, it is toric Morse on $(X_\Sigma, p^{-1}\{u_2\cdots u_m=0\})$ in the following sense.
\end{exa}
\begin{defin}\label{morsishtoric}
Let $V$ be a compact complex algebraic variety with a Weil divisor $D\supset \sing V$. A meromorphic function $f:V\DashedArrow\CP^1$ is said to be {\it toric Morse on} $(V,D)$, if it satisfies conditions (1-2) of Definition \ref{morsish} and the following instead of (3-4):

Near every point $x_0\in D$ 
the variety $V$ is locally isomorphic to an affine toric varitety $X$, so that under this isomorphism $D$ is a union of orbits, and 
one of the following takes place: 

3') if $f(x_0)\in\C^*$, then $X=X'\times\C$, and the map $f$ near $x_0$ factors through $X=X'\times\C\to \C$.

4') if $x_0\in\overline{\{f=0\}}$, then $\overline{\{f=0\}}$ is a closure of an orbit of $X$;

\end{defin}
\begin{rem}\label{remresol} In the setting of Example \ref{examerom}, we have by construction:

-- The function $\Phi$ is toric Morse and no indeterminacy points;

-- $p:\{\Phi=0\}\to\{\varphi=0\}$ is a homeomorphism.

In particular, verifying condition (3) 
amounts to checking that every non-zero level set 
$$\{\Phi=c\}=p^{-1}\{x_1^{d_1}-cx_2^{d_2}\cdots x_m^{d_m}=0\}\eqno{(*)}$$is transversal to every orbit of $X_\Sigma$. The latter is true because $\Sigma$ is the dual fan of the Newton polyhedron of the defining equation of $(*)$ (this explains how we invented our definition of $\Sigma$).
\end{rem}
We now define the resolution of the indeterminacy locus for an arbitrary Morse function $f$ on $(V,D)$. Near every point $z\in D\cap\overline{\{f=0\}}$, the function $f$ looks like the product of Example \ref{examerom} and $\C^{n-m}$ for some $m$. In the ring of germs of analytic functions on $(V,z)$, we define the monomial ideal $I_z$ generated by $x_1^{d_1}$ and $x_2^{d_2}\cdots x_m ^{d_m}$ (using the notation of Example \ref{examerom}). These ideals form a sheaf on $V$, and we define $p:\tilde V\to V$ as the resolution along this sheaf of ideals. 
\begin{utver}
In the setting of the preceding paragraph, we have the following:

-- the function $\tilde f:=f\circ p$ is toric Morse with no indeterminacy points;

-- $p:\{\tilde f=0\}\to\{f=0\}$ is a homeomorphism.
\end{utver}
\begin{proof}
Near every point $z\in D\cap\overline{\{f=0\}}$, this resolution is isomorphic to Example \ref{examerom}, so the statement follows by Remark \ref{remresol}.
\end{proof}
\subsection{Proofs}\label{ssproofsconn}
{\it Proof of Theorem \ref{toricmorse} and Remark \ref{realmilnor}.2.} By the preceding proposition, we can assume with no loss in generality that $f$ is toric Morse with no indeterminacy points. In this case condition (3') of Definition \ref{morsishtoric} implies:
$$
\mbox{Every } c\in\C^*\mbox{ is a regular value of }f\mbox{ restricted to every toric stratum of }D.\eqno{(*)}$$

This observation and Remark \ref{realmilnor}.1 allow application of  the standard Morse theoretic arguments to prove Theorem \ref{toricmorse}, except for Part 3 with $c=0$.

Let us prove the remaining part for small $c'$. This amounts to constructing a vector field $v$ on $U$ near ${\{f=0\}}$, retracting $\{|f|\leqslant c'\}$ to ${\{f=0\}}$. In a small neighborhood $U_z$ of a point $z\in \{f=0\}\setminus D$, such a vector field $v_z$ is known to exist. In a small neighborhood $U_z$ of a point $z\in \{f=0\}\cap D$, the variety $V$ is locally isomorphic to an affine toric variety with a dense torus $T$, and $\{f=0\}\cap U$ is its codimension 1 orbit $T'$. The sought vector field $v_z$ in $U_z$ is then the unique (up to scaling) $T$-invariant vector field on $T$ whose trajectories tend to $T'$. Now when we have the sought vector field $v_z$ in a small neighborhood $U_z$ of every point $z\in\{f=0\}$, we can use compactness of $V$ to choose a finite covering $\bigcup_i U_{z_i}\supset\{f=0\}$, and define the global vector field $v=\sum\psi_iv_{z_i}$ by means of a partition of unity $\psi_i:U_{z_i}\to\R,\,\sum_i\psi_i=1$.

Finally, let us prove Remark \ref{realmilnor}.2: if $W:=U\setminus\{f=0\}$ 
 is affine, then condition (2) of 
Definition \ref{morsishtoric} is redundant. Assume, towards, the contradiction, that the other conditions are satisfied, but the set of critical points in $W$ has a component of positive dimension. This component $C$ belongs to $\{f=c\}$ for some $c\in\C^*$. Since $W$ is affine, $C$ has a limit point $z\in D$. Since $z$ is not in the indetemrinacy locus of $f$ (which is empty), it is a critical point of $f$ on the stratum of $D$ containing $z$, which contradicts the observation $(*)$. $\hfill\square$

\vspace{1.5ex}

{\it Proof of Theorem \ref{th0connect}.}
i) The complete intersection $H_k$ has at most ICIS singularities.

\begin{proof}
By induction on $k$, the complete intersection $H_{k-1}$ has at most isolated singularities. If, towards the contradiction, the singular locus of the divisor $H_k\subset H_{k-1}$ has positive dimension (including the points where the divisor is not reduced), then it has a limit point $z\notin T$ in a tropical compactification. Near this point, $H_k$ cannot be smooth and reduced, contradicting Definition \ref{defschon0} of sch\"on.
\end{proof}

ii) $\beta_j(H_k)=\beta_j(T)$ for $j<n-k$.

\begin{proof}
Choose a smooth tropical compactification $X$ and a rational point $a$ in the interior of $N_k$ and positive $m\in\Z$ such that $am$ is a lattice point. The statement follows by Theorem \ref{toricmorse}, amended with Remarks \ref{realmilnor}.2 and \ref{remmorseicis} and applied to $f=f_k^m/x^{am}$, $V=\bar H_{k-1}\subset X$ with $D=V\setminus H_{k-1}$. By construction, $f$ is Morse on $(V,D)$:

-- we have $\{f=\infty\}\subset D$ because $a$ is in the interior of $N_k$, 

-- the set $\overline{\{f=0\}}$ is smooth near $D$ and transversal to it by the definition of sch\"on,

-- the variety $V$ has at most ICIS singularities by (i).

Thus Theorem \ref{toricmorse} applies and gives $\beta_j(H_k)=\beta_j(H_{k-1})$, which equals $\beta_j(T)$ by induction on $k$.
\end{proof}

iii) $H_k$ is irreducible unless it is a singular curve.

\begin{proof} 
Assume towards the contradiction that $H_k$ has at least two components $P$ and $Q$. Since $H_k$ is connected by (ii), there is a point $z\in P\cap Q$. Since it is a singular point of $H_k$, it is ICIS by (i). This contradicts with the fact that an ICIS of dimension greater than 1 is locally irreducible (e.g. because its link is connected). 
\end{proof}

{\it An excursus: a Morse theoretic proof of Example \ref{exascibkk}.1 and Theorem \ref{thkhconnect} for $k=1$.} When Bernstein proved formula (BKK) of Example \ref{exascibkk}.1 for $k=n$, he then deduced this formula and Theorem \ref{thkhconnect} for $k=1$ using Morse theory (according to \cite{dkh}, Remark 3.10), but never published this approach. We sketch what might be a reconstruction of his proof: the appeal to Morse theory is hidden in the reference to Theorem \ref{toricmorse}.

For $n>1$ and $k=1$ in the setting of Theorem \ref{thkhconnect}, assume with no loss of generality that the convex hull $P$ of $A_1\subset\Z^n$ is a full dimensional polytope, take a rational point $a$ in its interior, and apply Theorem \ref{toricmorse} in the following setting: $V$ is a smooth toric variety whose fan is compatible with $P$, $D$ is the complement to $\CC^n$, and the function $f$ is $f_1^m/x^{ma}$ for a nondegenerate Laurent polynomial $f_1\in\C^{A_1}$ and $m\in\Z$ such that $am\in\Z^n$.

Since $a$ is in the interior of $P$, we have $\overline{\{f=\infty\}}=D$. 
Since $f_1$ is nondegenerate, the closure $\overline{\{f=0\}}$ is smooth and transversal to every orbit of $V$, thus Theorem \ref{toricmorse} is applicable. 

First, we conclude that $\{f_1=0\}=\{f=0\}\subset\CC^n$ has the same Betti numbers as $\CC^n$ in all dimensions except for $n$ and $n-1$: in particular it has one connected component.

Second, we conclude that the Euler characteristics $(-1)^{n-1}\chi\{f_1=0\}=(-1)^{n-1}\chi\{f=0\}$ equals the sum of $\chi\CC^n=0$ and the total multiplicity of critical points of $f$ outside $\{f=0\}$. The latter equals the lattice volume of $P$, applying formula (BKK) to the square system of equations $x_i\partial f_1/\partial x_i-a_if_1=0,\,i=1,\ldots,n$, all supported at $A_1$. It is applicable to this system of equations, because the system is nondegenerate at infinity: no face $\Gamma\subsetneq P$ has $a$ in its affine span, so the restrictions $f_1|_\Gamma$ and $(\partial f_1/\partial x_i)|_\Gamma$ are linear combinations of the restrictions $(\partial f_1/\partial x_i-a_if_1)|_\Gamma$, so any common root of the latters would be a common root of the formers, contradicting the nondegeneracy of $f_1$.

We have proved formula (BKK) of Example \ref{exascibkk}.1 and Theorem \ref{thkhconnect} for $k=1$. 
\begin{rem}
Further, one can deduce formula (BKK) for $1<k<n$ from the case $k=1$: notice that the projection of the hypersurface $H=\{\sum_{i=1}^k \lambda_if_i(x)=0\}\subset\CP^{k-1}_\lambda\times\CC^n_x$ to $\CC^n_x$ has fibers equal to $\CP^{k-1}$ over the points of $\{f_1=\dots=f_k=0\}$, and $\CP^{k-2}$ over other points. Thus $e H$ equals the sought $e\{f_1=\dots=f_k=0\}$ by additivity of the Euler characteristics, and up to a sign equals the volume of the Newton polytope $\Delta$ of $\sum_{i=1}^k \lambda_if_i(x)$ (this is the case $k=1$ of the BKK formula which we have already proved). It remains to notice that the volume of $\Delta$ equals the sum of the mixed volumes of $A_i$'s in (BKK) for $k>1$.
\end{rem}

\section{Euler characteristics recursion}\label{sec5}

We assign the following to an SCI $H_k\subset\cdots\subset H_1$ in a torus $T$ with a cheracter lattice $M$:

\vspace{1ex}

-- its {\it defect} $\mu(H_k):=
(-1)^{n-k}(\mbox{the sum of the Milnor numbers of its singularities})$;

-- its smooth {\it tropical compactification} $X_\Sigma$ defined by a fan $\Sigma$ in $M^*$, as in Definition \ref{deftropcomp};

-- its {\it restriction} $(H_k)^l:=\bar H_k\cap O$, where $l\in M$ is the generator of a 1-dimensional cone in the fan $\Sigma$, and $O$ is the respective codimension 1 orbit of the toric variety $X_\Sigma$;

-- its {\it boundary data} $(l_{k,\alpha},m_{k,\alpha})$, where $\alpha$ is a component of $(H_{k-1})^{l_{k,\alpha}}$, and $m_{k,\alpha}$ is the multiplicity of $\alpha$ in the divisor of zeroes and poles of $f_k$, as in Definition \ref{defnewtpoyh};

-- its property of being Newtonian,
as in Definition \ref{defnewtsci};

-- in the Newtonian case, its {\it tropicaliztion} $(F_i,m_i)_{i=1,\ldots,k}$ as in Definition \ref{tropeq}.

\begin{rem}\label{remnewtind}
1. The restriction $(H_k)^l$ is a smooth SCI, and its boundary data can be combinatorially extracted from that of $H_k$ (c.f. Remark \ref{tropdescend});

2. $(H_k)^l$ is Newtonian if $H_k$ is Newtonian; $H_k$ is Newtonian if all $(H_i)^l$ are connected;

3. In particular, we can apply Theorem \ref{th0connect} to confirm connectivity of the SCIs $(H_k)^l$, and ultimately the Newtonian property of $H_k$, using only the boundary data of $H_k$ as the input.
\end{rem}

\subsection{The recursive theorem} 
We shall now express the Euler characteristics of SCI $H_k$ in terms of its 
tropicalization and the Euler characteristics of the restrictions $(H_i)^l$. 

Since the tropicalizations of the latters can be extracted from that of $H_k$ (Remark \ref{tropdescend}), this will allow us to compute the Euler characteristics of a SCI recursively by induction on the dimension of the ambient torus, using only its Newton data as the input. 

\begin{rem}\label{tropdescend}
While the recursion may look involved at a glance, it converges to simple explicit formulas like Theorem \ref{th0tropeuler}, or Theorems 1.4, 1.8 and 5.6 
in \cite{crit}. This is because the tropicalization of $(H_k)^l$ can be as follows extracted from that of Newtonian SCI $H_k$. 

Assume wlog that $m_k(l)=0$: we can achieve this, dividing the defining equation $f_k$ of the set $H_k$ by a suitable character of the torus (this operation does not change $H_k$). Then there exists the unique tropical complete intersection $(F'_i,m'_i)$ in the quotient lattice $M^*/l$, whose lift to $M^*$ coincides with the tropicalization $(F_i,m_i)$ of $H_k$ in a small neighborhood of $l$. This $(F'_i,m'_i)$ is the tropicalization of $(H_i)^l$. 
\end{rem}

\begin{theor}\label{th0euler}
For a Newtonian SCI $C_m\subset\cdots\subset C_1$, extend its defining function $f_i$ on $C_{i-1},\, i\leqslant m$, to a Laurent polynomial $F_i:T\to\C$. 
Let $P_i\subset M$ be a lattice polytope containing the monomials of $F_i$ (e.g. the Newton polytope of $F_i$ works well for most applications), and $\tilde F_i$ be a generic Laurent polynomial with the same Newton polytope $P_i$. 

Denote the tuple $(F_1,\ldots,F_{i-1},\tilde F_{i+1},\ldots,\tilde F_m)$ by $F^i$ and the support function of $P_i$ by $\tilde m_i$. Choose a smooth tropical compactification $X_\Sigma$ of $C_m$, such that $\tilde m_i$ is linear on every cone of $\Sigma$. 
Then the Euler characteristics $e(C_m)$ equals $e(F^0=0)+\mu(C_m)$ minus
$$\sum_{i=1}^m\sum_l (\tilde m_i(l)-m_i(l))\cdot\Bigl(e(F^i=0)^l-e(F^i=F_i=0)^l-e(F^i=\tilde F_i=0)^l+e(F^i=F_i=\tilde F_i=0)^l\Bigr),$$
where $l$ runs over the generators of all 1-dimensional cones in the fan $\Sigma$.

\end{theor}
\begin{rem}
1. This can be regarded as a recursive way to compute $e(C_m)$, because the other terms in the equality are known in the following sense: 

-- $e(F^0=0)$ can be computed in terms of $P_i$'s by formula (BKK) of Example \ref{exascibkk}.1 (which is itself a corollary of Theorem \ref{th0euler}, see Section \ref{sproof2}); 

-- the four terms $e(\cdots)^l$ can be found from the tropicalizations of $C_m$, applying the same Theorem \ref{th0euler} by induction on dimension $n$ (this application is possible thanks to Remark \ref{tropdescend}).

2. Theorem \ref{th0tropeuler}, giving a more convoluted formula for the Euler characteristics of a SCI, will be deduced from this one. However, this recursive algorithm is more flexible than the convoluted formula: for instance, it gives a simpler way to prove Theorem 1.8 
in \cite{crit} on Euler obstructions, because the definition of the Euler obstruction is recursive as well.

3. A sufficient genericity assumption for $\tilde F_i$'s in Theorem \ref{th0euler} is that the following are SCIs:
$$F_1=\cdots=F_{i-1}=\tilde F_{i+1}=\cdots=\tilde F_m=0,\, i=0,\ldots,m-1;$$
$$F_1=\cdots=F_{i}=\tilde F_{i+1}=\cdots=\tilde F_m=0,\, i=1,\ldots,m-1;$$
$$F_1=\cdots=F_{i}=\tilde F_{i}=\cdots=\tilde F_m=0,\, i=1,\ldots,m.$$
\end{rem}
The proof of the theorem (making do with this particular genericity assumption for $\tilde F_i$'s) will follow the original deformation strategy that was first employed in \cite{bernst} to count points and then in \cite{varch} to count Euler characteristics. 

What seems to be new in Theorem \ref{th0euler} is not a novel approach, but an observation that the classical approach allows induction on the dimension in a much higher generality.

\subsection{Proof of Theorem \ref{th0euler}} The sought equality is the sum over all $i$ of the following ones. 
\begin{lemma}[The recursive lemma] \label{l0euler}For every $i=1,\ldots,m$, we have
$$e(F_1=\cdots=F_{i-1}=\tilde F_{i}=\cdots=\tilde F_m=0)-e(F_1=\cdots=F_{i}=\tilde F_{i+1}=\cdots=\tilde F_m=0)=$$ $$=\sum_l (\tilde m_i(l)-m_i(l))\cdot\Bigl(e(F^i=0)^l-e(F^i=F_i=0)^l-e(F^i=\tilde F_i=0)^l+e(F^i=F_i=\tilde F_i=0)^l\Bigr),$$ with $\mu(C_m)$ subtracted in case $i=m$. 
\end{lemma}
We precede the proof with a remark and an example.
\begin{rem}\label{rempermut}
1. If $g_1=\cdots=g_m=0$ is a SCI, and $\sigma$ is a permutation of $\{1,\ldots,m\}$, then, in general, $g_{\sigma(1)}=\cdots=g_{\sigma(m)}=0$ need not be a SCI. Furthermore, even if both are SCIs, still the boundary data of $g_i$ in these two SCIs differs uncontrollably (i.e. the difference is not determined by the boundary data of all of the $g_1,\ldots,g_m$ in the first or the second SCI).

However, the boundary data of $g_j$ is the same in both SCIs, if $\sigma(j)=j$ and $i<j\Rightarrow\sigma(i)<\sigma(j)$, or if $g_j$ is a generic polynomial with a given Newton polytope.

Moreover, if $\sigma$ is a transposition $(j,j+1)$, and $g_j$ is a generic polynomial with a given Newton polytope, then the boundary data of $g_{j+1}$ is as well the same in both SCIs.

2. In particular, for the SCI in our lemma, the boundary data of every equation in $F_1=\cdots=F_{i}=\tilde F_{i+1}=\cdots=\tilde F_m=0$ is the same as in $F^i=F_i=0$ (where $F_i$ is sent to the very end, and the order of the rest of the equations is unchanged). This is because all of the the equations located between the old and the new position of $F_i$ are generic polynomials with prescribed Newton polytopes. We shall use this seamlessly in the proof of the lemma.
\end{rem}
For germs of analytic functions $g$ and $\tilde g:(\C^q,z)\to(\C,0)$ and a Zariski open set $U\subset\C^q$, the {\it Milnor fiber at $z$ of the meromorphic function} $g/\tilde g|_U$ is the intersection of $U$, a $\delta$-neighborood of $z\in\C^q$ and the hypersurface $\{g=\varepsilon\cdot\tilde g\}$ for $\varepsilon\ll\delta\ll 1$. Its topological type does not depend on $\varepsilon$ and $\delta$, so we denote it by $\mu_z(g/\tilde g|_U)$.
\begin{exa}\label{examilntori}
The Milnor fiber $\mu_0(x_1^{m_1}\cdots x_q^{m_q}|_{x_1\cdots x_k\ne 0})$ for $m_\bullet\geqslant 0$ is homotopy equivalent to the disjoint union of $GCD(m_1,\ldots,m_q)$ tori of dimension $\Bigl|[1,\ldots,k]\cup\{j\,|\,m_j>0\}\Bigr|-1$.

The Milnor fiber $\mu_0(x_2^{m_2}\cdots x_q^{m_q}/x_1|_{x_2\cdots x_k\ne 0})$ for $m_\bullet\geqslant 0$ is homotopy equivalent to a torus of dimension $k-1$ (where negative dimension means empty set).
\end{exa}

{\it Proof of Lemma \ref{l0euler}.} The sought difference $e(F^i=\tilde F_i=0)-e(F^i= F_i=0)$ equals $e(F^i=F_i-\varepsilon\tilde F_i=0)-e(F^i= F_i=0)$ for small $\varepsilon$, because both $\tilde F_i$ and $F_i-\varepsilon\tilde F_i$ are generic polynomials with the Newton polytope $P_i$ under the assumptions of this lemma.

The latter difference equals the Euler characteristics integral of $e(\mu_z(F_i/\tilde F_i|_{C_{i-1}}))$ over $z\in\bar V_{i-1}\setminus C_{i-1}$ (minus $\mu(C_m)$ in case $i=m$).
Every such Milnor fiber has the form of Example \ref{examilntori} (in a suitable coordinate system near $z$), so it is a disjoint union of tori. These tori have positive dimension unless, for a generator $l$ of some ray in the fan $\Sigma$, the point $z$ belongs to $$Z_l:=\{F^i=0\}^l\setminus(\{F^i=F_i=0\}^l\cup\{F^i=\tilde F_i=0\}^l).$$ 
In the latter case, the number of 0-dimensional tori is $\tilde m_i(l)-m_i(l)$.

Thus the sought integral of $e(\mu_z(F_i/\tilde F_i|_{C_{i-1}}))$ over $z\in\bar V_{i-1}\setminus C_{i-1}$ equals $\sum_l (\tilde m_i(l)-m_i(l))\cdot e(Z_l).$ This coincides with the statement of lemma by additivity of Euler characteristics. \hfill$\square$

\begin{rem}\label{remnonnewt}
Assume we wish to compute the Euler characteristic of an SCI which is not Newtonian, extending Lemma \ref{l0euler}, Theorem \ref{th0euler} and ultimately Theorem \ref{th0tropeuler} to arbitrary SCIs.
In the proof of the extended Lemma \ref{l0euler}, we shall need Euler characteristics of the individual components of $\{F^i=0\}^l$ rather than their sum, because the Milnor number $e(\mu_z(F_i/\tilde F_i|_{C_{i-1}}))$ as a function of $z$ may take distinct values on different components of $\{F^i=0\}^l$.

So, in order to make the induction's ends meet, the extended Lemma \ref{l0euler} should compute the difference of the Euler characteristics of the respective components of the SCIs $F^i=F_i=0$ and $F^i=F_i-\varepsilon\tilde F_i=0$, but for this we need a correspondence between these components: the components of the latter should not split into more components as $\varepsilon$ tends to 0.

We can ensure this easily, given the positive answer to Question \ref{conjflat}: the components of the SCI $F^i=F_i=0$ then belong to distinct cosets of a subtorus $T_i$, so let $M_i\subset M$ be the sublattice of characters that equal 1 on $T_i$, and choose $P_i$ to be a sufficiently large polytope in $M_i$.

\end{rem}

\subsection{Deducing Theorem \ref{th0tropeuler} from Theorem \ref{th0euler}} Part 1 follows from Proposition \ref{schontopol}. 
For Part 2, see e.g. \cite{mikhicm}, \cite{rau07}, \cite{k09}, \cite{mmjpoly} for background on tropical intersection theory.

\begin{theor}\label{th1tropeuler}
In the setting of Theorem \ref{th0euler}, we have $$e(F_1=\cdots=F_{i}=\tilde F_{i+1}=\cdots=\tilde F_m=0)=\frac{\delta\tilde m_k}{1+\delta\tilde m_k}\cdots\frac{\delta\tilde m_{i+1}}{1+\delta\tilde m_{i+1}}\frac{\delta m_{i}}{1+\delta m_{i}}\cdots\frac{\delta m_{1}}{1+\delta m_{1}}$$
with $\mu(C_m)$ added in case $i=m$.
\end{theor}
\begin{proof} 
Note that the statement for smaller $i$ is the special case of the same statement for larger $i$. We formulate it so, because the proof will proceed by induction on $i$, assuming that it is already known for smaller $i$ and for Newtonian SCIs in tori of smaller dimension (the base of induction on both parameters is given by Example \ref{exaropnewt}).

By this inductive assumption, the sought equality is equivalent to the first one in the chain
$$e(F^i=\tilde F_i=0)-e(F^i= F_i=0)=$$
$$=\frac{\delta\tilde m_k}{1+\delta\tilde m_k}\cdots\frac{\delta\tilde m_{i}}{1+\delta\tilde m_{i}}\frac{\delta m_{i-1}}{1+\delta m_{i-1}}\cdots\frac{\delta m_{1}}{1+\delta m_{1}}-\frac{\delta\tilde m_k}{1+\delta\tilde m_k}\cdots\frac{\delta\tilde m_{i+1}}{1+\delta\tilde m_{i+1}}\frac{\delta m_{i}}{1+\delta m_{i}}\cdots\frac{\delta m_{1}}{1+\delta m_{1}}=$$ 
$$=\frac{\delta\tilde m_k}{1+\delta\tilde m_k}\cdots\frac{\delta\tilde m_{i+1}}{1+\delta\tilde m_{i+1}}\left(\frac{\delta\tilde m_i}{1+\delta\tilde m_{i}}-\frac{\delta m_i}{1+\delta m_{i}}\right)\frac{\delta m_{i-1}}{1+\delta m_{i-1}}\cdots\frac{\delta m_{1}}{1+\delta m_{1}}=$$
$$=\frac{\delta\tilde m_k}{1+\delta\tilde m_k}\cdots\frac{\delta\tilde m_{i+1}}{1+\delta\tilde m_{i+1}}\frac{\delta(\tilde m_i-m_i)}{(1+\delta\tilde m_{i})(1+\delta m_{i})}\frac{\delta m_{i-1}}{1+\delta m_{i-1}}\cdots\frac{\delta m_{1}}{1+\delta m_{1}}.\eqno{(*)}$$

To prove it, let us consider a 1-dimensional tropical fan $R$, whose rays are 1-dimensional cones $\langle l\rangle\in\Sigma$, and whose weights are given by the huge bracket on the right hand side of Lemma \ref{l0euler} (incidentally, it is the 1-dimensional tropical characteristic class of the set $(F^i=0)\setminus(F_i\cdot\tilde F_i=0)$ in the sense of \cite{e13}, though we do not use this in the proof). We can compute this huge bracket on the right of Lemma \ref{l0euler} by induction on the dimension of the ambient torus, because we can use Remark \ref{tropdescend} to extract the tropicalizations of the participating SCIs from $\tilde m_\bullet$ and $m_\bullet$. The result is as follows.

$$R=\frac{\delta\tilde m_k}{1+\delta\tilde m_k}\cdots\frac{\delta\tilde m_{i+1}}{1+\delta\tilde m_{i+1}}\frac{\delta m_{i-1}}{1+\delta m_{i-1}}\cdots\frac{\delta m_{1}}{1+\delta m_{1}}-\frac{\delta\tilde m_k}{1+\delta\tilde m_k}\cdots\frac{\delta\tilde m_{i}}{1+\delta\tilde m_{i}}\frac{\delta m_{i-1}}{1+\delta m_{i-1}}\cdots\frac{\delta m_{1}}{1+\delta m_{1}}-$$
$$-\frac{\delta\tilde m_k}{1+\delta\tilde m_k}\cdots\frac{\delta\tilde m_{i+1}}{1+\delta\tilde m_{i+1}}\frac{\delta m_{i}}{1+\delta m_{i}}\cdots\frac{\delta m_{1}}{1+\delta m_{1}}+\frac{\delta\tilde m_k}{1+\delta\tilde m_k}\cdots\frac{\delta\tilde m_{i}}{1+\delta\tilde m_{i}}\frac{\delta m_{i}}{1+\delta m_{i}}\cdots\frac{\delta m_{1}}{1+\delta m_{1}}=$$
$$=\frac{\delta\tilde m_k}{1+\delta\tilde m_k}\cdots\frac{\delta\tilde m_{i+1}}{1+\delta\tilde m_{i+1}}\left(1-\frac{\delta\tilde m_{i}}{1+\delta\tilde m_{i}}-\frac{\delta m_{i}}{1+\delta m_{i}}+\frac{\delta\tilde m_{i}}{1+\delta\tilde m_{i}}\frac{\delta m_{i}}{1+\delta m_{i}}\right)\frac{\delta m_{i-1}}{1+\delta m_{i-1}}\cdots\frac{\delta m_{1}}{1+\delta m_{1}}=$$
$$=\frac{\delta\tilde m_k}{1+\delta\tilde m_k}\cdot\frac{\delta\tilde m_{i+1}}{1+\delta\tilde m_{i+1}}\frac{1}{(1+\delta\tilde m_{i})(1+\delta m_{i})}\frac{\delta m_{i-1}}{1+\delta m_{i-1}}\cdots\frac{\delta m_{1}}{1+\delta m_{1}}.$$
Thus the whole right hand side of Lemma \ref{l0euler} is the weight of the 0-dimensional tropical fan
$$\delta(\tilde m_i-m_i)\cdot R=\delta(\tilde m_i-m_i)\frac{\delta\tilde m_k}{1+\delta\tilde m_k}\cdot\frac{\delta\tilde m_{i+1}}{1+\delta\tilde m_{i+1}}\frac{1}{(1+\delta\tilde m_{i})(1+\delta m_{i})}\frac{\delta m_{i-1}}{1+\delta m_{i-1}}\cdots\frac{\delta m_{1}}{1+\delta m_{1}},$$
which is the right hand side of the sought equality $(*)$. Thus Lemma \ref{l0euler} proves the sought equality $(*)$. Theorem \ref{th1tropeuler} and thus \ref{th0tropeuler} is proved. \end{proof}

\subsection{An excursus: deducing formula (BKK) of Example \ref{exascibkk}.1 from Theorem \ref{th0euler}} \label{sproof2}
By Theorem \ref{bkktopol}, it is enough to prove the equality for any particular nondegenerate system of equations $g=0$. We shall do this in the case when $g=0$ is SCI (which is a stronger assumption).

If a codimension 1 orbit $O$ of a tropical compactification $X_\Sigma$ corresponds to the 1-dimensional cone of the fan $\Sigma$, generated by primitive $\gamma$, then denote $A_i^\gamma$ by $A_i^O$. Assume with no loss in generality that $0\in A_1$. Apply Theorem \ref{th0euler} to generic $F_i$'s with the Newton polytopes $P_i:=\conv A_i$ for $i=1,\ldots,m$, $m:=k$, $F_1:=1$, and $\tilde F_i:=g_i$ for $i=2,\ldots,m$. Since $\{F_1=0\}=\varnothing$, we have $e(C_m)=0$. Writing $e_{n-1}(A_1^O,\ldots,A_m^O)$ instead of $e(\tilde F_1=\cdots=\tilde F_m=0)^O$ and $e_{n-1}(\{0\}, A_2^O,\ldots,A_m^O)$ instead of $e(F_1=\cdots=0)^O$ by induction on the dimension $n$ for every term in the sum, the sought Euler characteristics $e(F^0=0)$ equals 
$$\sum_O\tilde m_1^O\cdot\bigl(e_{n-1}(A_2^O,\ldots,A_m^O)-e_{n-1}(A_1^O,\ldots,A_m^O)\bigr)=$$ $$=(-1)^{n-k}\sum_\gamma \max\gamma|_{A_1}\cdot \sum_{d_1+\ldots+d_k=n\atop d_1,\ldots,d_k>0} (A_1^\gamma)^{d_1-1}(A_2^\gamma)^{d_2}\cdots (A_k^\gamma)^{d_k}=(-1)^{n-k}\sum_{d_1+\ldots+d_k=n\atop d_1,\ldots,d_k>0} A_1^{d_1}\cdots A_k^{d_k}.$$
Here, in the first equality, we pass from summation over codimension 1 orbits of $X_\Sigma$ (i.e. over primitive $\gamma\in(\Z^n)^*$ generating 1-dimensional cones of $\Sigma$) to summation over all primitive $\gamma\in(\Z^n)^*$. This is valid, because, once $\gamma$ belongs to the relative interior of a $q$-dimensional cone $C\in\Sigma$ for $q>1$, all $A_i^\gamma$'s can be shifted to the codimension $q$ normal plane of $C$, and thus the corresponding term of the sum vanishes. The last equality holds term by term, because the left hand side term for given $d_1,\ldots,d_k$ is the ``base$\times$height'' formula for the mixed volume $A_1^{d_1}\cdots A_k^{d_k}$ (see e.g. \cite{ew}, Theorem 4.10).

\section{Discussion}\label{sdiscus}

\begin{exa}\label{exaempty}
For generic Laurent polynomials $f, g\in\C^{[-a,a]}$, consider a Newtonian SCI $$f(x)g(y)-1=f(x)(y^{-p}+y^q)=0.$$ 
The tropicalization of the second equation is not convex, and its Newton polyhedron is empty whenever $p+q<2a$, because the boundary data of the second equation includes pairs $((0,-1),p-a)$ and $((0,1),q-a)$. 

In particular, the number of solutions $2(p+q)a$ of this system is not related to the mixed area of their Newton polyhedra $[-a,a]^2$ and $\varnothing$ or Newton polytopes $[-a,a]^2$ and $[-a,a]\times[-p,q]$. However, Theorem \ref{th0tropeuler} works well: the tropicalization of the first equation $m_1$ is the support function of $[-a,a]^2$, so $\delta m_1$ consists of the rays generated by $(0,\pm 1)$ and $(\pm 1, 0 )$ with multiplicities $2a$. The tropicalization of the second equation $m_2$ equals $a,a,q-a,p-a$ on these four vectors, so $\delta (m_2\cdot\delta m_1)=\frac{\delta^2}{2!}(m_1m_2)$ has multiplicity $2a\cdot (a+a+(q-a)+(p-a))=2(p+q)a$.

Of course, the same complete intersection set can be defined by a nondegenerate equations $f(x)g(y)-1=(y^{-p}+y^q)=0$, so the number of solutions is just the mixed area of their Newton polygons $[-a,a]^2$ and $\{0\}\times [p,q]$.

\end{exa}
With this example in mind, we highlight some natural open questions on SCIs. 

\subsection{Betti numbers of SCIs} 
What happens to Theorem \ref{th0connect} if the Newton polyhedra of a SCI are not of full dimension? In the classical case of nondegenerate comple intersections, Theorem \ref{thkhconnect} gives a complete answer, after which we model our open questions.
\begin{quest}\label{qqqessent0}
Given a SCI $V\subset T$ with the Newton polyhedra $N_1,\ldots,N_k$, let $d$ be the minimum of numbers $\dim(N_{i_1}+\cdots+N_{i_q})-q$ over all tuples $i_1<\cdots<i_q$. Is it true that the Betti numbers of dimension smaller than $d$ coincide for $V$ and $T$?
\end{quest}
\begin{rem}
One can try to deal with it using the following significant extension of Theorem \ref{toricmorse}.1: if the critical points of $f$ at the level set $\{\varphi=c\}$ form a set of positive complex dimension $d$, then the set $\{\varphi\leqslant c\}\cap U$ is homotopy equivalent to $\{\varphi<c\}\cap U$ with several cells of real dimension $\dim_\C U-d$ or higher attached.
\end{rem}
Anyway, 
there is no chance to obtain a general exact criterion for connectivity of a SCI in terms of its Newton polyhedra: in Example \ref{exaempty}, the number of components (i.e. points) of the SCI depends on $p$ and $q$, while the Newton polyhedra do not remember these numbers. 

At best, we might expect such a criterion for a Newtonian SCI, in terms of the tropicalizations of its defining equations.
For a piecewise linear function $m$ on a tropical fan in $\R^n$, define its lineality space $L(m)\subset\R^n$ as the maximal space such that $m$ decomposes into 
(a lift of a piecewise linear function on $\R^n/L(m))+($a linear function on $\R^n)$. 
\begin{quest}\label{qqqessent}
1. Given a Newtonian SCI $V\subset T$ with the tropicalizations of the equations $m_1,\ldots,m_k$, let $d$ be the minimum of numbers $\dim(L(m_{i_1})+\cdots+L(m_{i_q}))-q$ over all tuples $i_1<\cdots<i_q$. Is it true that the Betti numbers of dimension smaller than $d$ coincide for $V$ and $T$, and in dimension $d$ we may have a strict inequality?

2. Are the Betti numbers of a SCI determined by its boundary data? (By the Betti number of dimension $i$ we mean the dimension of the $i$-th topolgical homology group modulo infinity.)
\end{quest}

\subsection{Geometry of tropical complete intersections} The rest of the questions are motivated by our main intention behind this paper: to construct subvarieties of the algebraic tori with prescribed topology, starting from a rich yet manageable combinatorial data.

We intend to use a tropical complete intersection $\frak{V}$ as such a combinatorial data, and to construct a Newtonian SCI $V$ with the tropicalization $\frak{V}$. Many geometric characteristics $E$ of a variety $V$ can be read off from $\frak{V}$, as
Theorems \ref{th0connect} and \ref{th0tropeuler} suggest: this means that we can define $E(\frak{V})$ for a tropical complete intersection $\frak{V}$ so that it equals $E(V)$, once $\frak{V}$ is the tropicalization of $V$.
We shall call such characteristics $E$ tropicalizable. (Recall that the fundamental group is an example of a non-tropicalizable topological invariant \cite{rybn}.)

Constructing an SCI with prescribed geometric characteristics $E$ which is tropicalizable (such as Euler characteristics) is done in two steps.

1. To find a tropical complete intersection $\frak{V}$ such that $E(\frak{V})$ takes the sought value.

2. To check that the space of Newtonian SCIs having the prescribed tropicalization $\frak{V}$ is non-empty, i.e. there exists an SCI with the sought tropicalization.

The two steps of this plan motivate the next two groups of questions: about topological invariants of tropical complete intersections, and about moduli spaces of SCIs.

Recall that a tropical complete intersection on a tropical fan $F_0$ in $\R^n$ is a sequence of nested fans $F_i$ of dimensions $\dim F_0-i$ and continuous piecewise linear functions $m_i:F_{i-1}\to\R$, whose corner loci equal $F_i$. 
Define its Euler characteristics by the formula $(**)$ of Theorem \ref{th0tropeuler}:
$$\frac{\delta m_k}{1+\delta m_k}\cdot\ldots\cdot\frac{\delta m_1}{1+\delta m_1}\cdot F_0.\eqno{(*)}$$
\begin{rem}
1. In the introduction, we restricted to the case $F_0=\R^n$.

2. If $k$ is the codimension of $F_0$, we call this quantity $(*)$ the {\it multiplicity} of the 0-dimensional tropical complete intersection.

3. Similarly to algebraic geometry over fields, it is natural to call a $k$-dimensional fan $F_0$ {\it smooth}, if there exist $m_1,\ldots,m_k$ for which the multiplicity equals 1. It is a priori wider than the usual notion of smoothness for tropical fans (which additionally assumes that $m_1=\ldots=m_k=($tropical linear function$)$, see \cite{fink}).

\end{rem}
\begin{quest}
1. 
Every $d$-dimensional local complete intersection in the algebraic torus has the Euler characteristics of the expected sign $(-1)^d$ or equals 0 (\cite{huh}, \cite{bw}). How to combinatorially define a subclass $\mathcal{P}$ of tropical complete intersections, which includes tropicalizations of Newtonian SCIs, and such that the Euler characteristics of every $d$-dimensional tropical complete intersection from $\mathcal{P}$ has sign $(-1)^d$ or equals 0?

2. How to classify tropical complete intersections in $\mathcal{P}$ with a small absolute value of Euler characteristics (such as 0, 1 or 2)? In particular, how to classify smooth tropical fans?

3. Is the classification finite for any given value of the Euler characteristics, if we restrict to ``non-homogeneous'' tropical complete intersections (such that $\dim(L(m_{i_1})+\cdots+L(m_{i_q}))>q$ for all $i_1<\cdots<i_q$)?
\end{quest}
\begin{rem}
While we do not know good candidates for the class $\mathcal{P}$, questions 2\&3 make sense even for a placeholder class of all $(F_0,m_1,\ldots,m_k)$ such that all cones of the fan $F$ has positive weights, and all functions $m_i$ are support functions of lattice polytopes. It does not include tropicalizations of all Newtoinian SCIs, but already with this understanding of $\mathcal{P}$ our questions have motivation in the geometry of lattice polytopes and matroids.

\end{rem}

\subsection{Deformations of SCIs} 

\begin{quest}
Find a positive formula for the Euler characteristics of a Newtonian SCI $C$ in terms of its tropicalization (i.e. a formula in which every summand has the same sign $(-1)^{\dim C}$ as the result). The question is open even in the simplest special cases such as Theorems 1.4, 1.8 and 5.6 in \cite{crit}.
\end{quest}
For a nondegenerate complete intersection $C\subset\CC^n$ defined by equations $f_i\in\C^{A_i}$, there is a classical construction of toric degeneration, tracing back e.g. to Viro's patchworking: pick ``generic'' sets $B_i\subset \Z^n\times\Z$ projecting to $A_i\subset\Z^n$, non-degenerate polynomials $g_i\in\C^{B_i}$ on $\CC^n_x\times\CC_t$ such that $g_i(x,1)=f_i(x)$, and a toric variety $X_\Sigma$ whose fan is compatible with $B_i$'s. Then, as $t\to 0$, the family of complete intersections $C_t\subset\CC^n$ defined by the equations $g_i(x,t)=0$ degenerates to the union of complete intersections $C^l$ given by the equations $g_i^l(x,1)=0$, where $l$ runs over the primitive generators of the rays in the fan $\Sigma$.

One application is that an additive nonnegative characteristic $E(C)$ splits into simple nonnegative summands $E(C^l)$ (this is literally true if $E$ is the Euler characteristics, and adapts to many other settings).

\begin{quest} Given a toric variety $X\supset\CC^n_x\times\CC_t$ and its subvariety $Y:=\{t=0\}$,
find general sufficient conditions on a Newtonian SCI $f_1=\cdots=f_k=0$ with a tropical compactigfication $Y$ under which $f_i$'s extend to a Newtonian SCI $g_1=\cdots=g_k=0$ with a tropical compactigfication $X$. This would in particular be applicable to the preceding question.
\end{quest}

\subsection{Moduli of SCIs} 
Let the fans $F_0=\R^n,F_1,\ldots,F_k$ and piecewise linear functions $m_1,\ldots,m_k$ on them define a tropical complete intersection. Consider the set $M_{m_1,\ldots,m_k}$ of all Newtonian SCIs whose defining equations tropicalize to $m_1,\ldots,m_k$.

Note that fibers of the forgetful projection $M_{m_1,\ldots,m_k}\to M_{m_1,\ldots,m_{k-1}}$ have the natural structure of Zariski open subsets of vector spaces.
\begin{quest}\label{qqqmoduli}
1. Do the fibers have the same dimension?
If so, then how to express this dimension in terms of $(m_1,\ldots,m_k)$?

2. How to introduce the structure of an algebraic variety on $M_{m_1,\ldots,m_k}$?
Will it turn the forgetful projection into a vector bundle projection on a Zariski open subset?

The answer to these questions may be negative. An example of such a tropical complete intersection can be extracted from  \cite{rybn}: in case of positive answers, the moduli space would be connected, any two SCIs with a given tropicalization could be connected with a continuous deformation (i.e. a path in the moduli space), and thus diffeomorphic by Theorem \ref{th0tropeuler}.1.

3. The simplest and most important special case of these question is whether $M_{m_1,\ldots,m_k}$ is non-empty for any given tropical complete intersection $(F_i,m_i)$, i.e. whether there exists a Newtonian SCI with a prescribed tropicalization.

\end{quest}

\vspace{1cm}
\noindent
London Institute for Mathematical Sciences, UK \\
\textit{Email}: aes@lims.ac.uk

\vspace{1ex}
























\begin{thebibliography}{xxxxxxxx}

\bibitem[AR07]{rau07} L. Allermann, J. Rau, {\it First Steps in Tropical Intersection Theory}, Math Z 264 (2010) 633-670, arXiv:0709.3705

\bibitem[AAGGL15]{gcicy}
L. B. Anderson, F. Apruzzi, X. Gao, J. Gray, S.-J. Lee, {\it A new construction of Calabi-Yau manifolds: Generalized CICYs} Nucl. Phys. B906 (2016) 441–496 arXiv:1507.03235

\bibitem[BB94]{bb94}
V. Batyrev, L. Borisov, {\it On Calabi-Yau Complete Intersections in Toric Varieties}, in Higher dimensional complex varieties, 39–65, de Gruyter, Berlin, 1996, arXiv:alg-geom/9412017

\bibitem[BH22]{bh}
P. Berglund, T. H\"ubsch, {\it Hirzebruch Surfaces, Tyurin Degenerations and Toric Mirrors: Bridging Generalized Calabi-Yau Constructions}, arXiv:2205.12827

\bibitem[B75]{bernst} 
D. N. Bernstein, {\it The number of roots of a system of equations}, Functional Anal. Appl. 9 (1975) 183--185. 

\bibitem[BW14]{bw} 
N. Budur, B. Wang, {\it The Signed Euler Characteristic of Very Affine Varieties}, IMRN 2015 (2015) 5710–5714, arXiv:1403.4371.


\bibitem[CHKS04]{chks}
F. Catanese, S. Hoşten, A. Khetan, B. Sturmfels, {\it The maximum likelihood degree}, American Journal of Mathematics 128 (2006) 671-697, arXiv:math/0406533

\bibitem[CKP14]{ckp}
T. Coates, A. Kasprzyk, T. Prince, {\it Four-dimensional Fano toric complete intersections}, Proc. R. Soc. A 471 (2015) 20140704, arXiv:1409.5030.

\bibitem[DKh87]{dkh}
V. Danilov, A. Khovanskii, {\it Newton polyhedra and an algorithm for computing Hodge–Deligne numbers}, Math. USSR-Izv 29:2 (1987) 279-298.

\bibitem[DP01]{dp}
A. Dimca, S. Papadima, {\it Hypersurface complements, Milnor fibers and higher homotopy groups of arrangements}, Ann. of Math. 158 (2003) 473-507, arXiv:math/0101246

\bibitem[DHOST13]{dhost}
J. Draisma, E. Horobeţ, G. Ottaviani, B. Sturmfels, R. R. Thomas, {\it The Euclidean distance degree of an algebraic variety}, Foundations of Comput. Math., 16 (2016) 99-149, 	arXiv:1309.0049









\bibitem[E10]{mmjpoly} A. Esterov, {\it Tropical varieties with polynomial weights and corner loci of piecewise polynomials}, Mosc. Math. J., 12:1 (2012), 55-76, arXiv:1012.5800

\bibitem[E11]{adv} A. Esterov, {\it The discriminant of a system of equations}, Adv. Math. 245 (2013) 534--572, arXiv:1110.4060

\bibitem[E13]{e13}A. Esterov, {\it Characteristic classes of affine varieties and Plucker formulas for affine morphisms}, Journal of the EMS, 20 (2018) 15-59, arXiv:1305.3234

\bibitem[EL22]{symm} A. Esterov, L. Lang, {\it Bernstein-Kouchnirenko-Khovanskii with a symmetry} arXiv:2207.03923

\bibitem[E24e]{crit} A. Esterov,
{\it Slightly degenerate Bernstein--Kouchnirenko-Khovanskii} arXiv:????.???? (pending the number)

\bibitem[E]{ew}
G. Ewald, {\it Combinatorial convexity and algebraic geometry}, GTM 168, Springer (1996). 












\bibitem[F10]{fink}
A. Fink, {\it Tropical cycles and Chow polytopes}, Beitr\"age zur Algebra und Geometrie 54 (2013) 13-40, arXiv:1001.4784

\bibitem[FS94]{fs}
W. fulton, B. Sturmfels, {\it Intersection theory on toric varieties}, Topology
36 (1997) 335-353, arXiv:alg-geom/9403002






\bibitem[H12]{huh}
J. Huh, {\it The maximum likelihood degree of a very affine variety}, Compos. Math. 149 (2013) 1245–1266, arXiv:1207.0553.

\bibitem[HS13]{hs}
J. Huh, B. Sturmfels, {\it Likeligood geometry}, in Combinatorial Algebraic Geometry, Lecture Notes in Mathematics 2108 (2014), arXiv:1305.7462

\bibitem[K09]{k09}
E. Katz {\it Tropical intersection theory from toric varieties}, Collectanea Mathematica 63 (2012) 29-44, 	arXiv:0907.2488

\bibitem[KP07]{kppoly}  E. Katz, S. Payne, {\it Piecewise polynomials, Minkowski weights, and localization
on toric varieties}, Algebra and Number Theory Journal, 2 (2008) 135-155, 	arXiv:math/0703672


\bibitem[K78]{kh75} 
A. G. Khovanskii, {\it Newton polyhedra and the genus of complete intersections}, Functional Anal. Appl., 12 (1978) 38--46


\bibitem[K77]{khovcomp}  
A. G. Khovanskii, {\it Newton polyhedra and toroidal varieties}, Functional Analysis and Its Applications, 11 (1977) 289--296.

\bibitem[K16]{kh15}  
A. G. Khovanskii, {\it Newton polytopes and irreducible components of complete intersections}, Izvestiya: Mathematics, 80 (2016) 263--284

\bibitem[L]{lij}
E. J. N. Looijenga, {\it Isolated Singular Points on Complete Intersections}, LMS Lecture Notes 77

\bibitem[MS]{ms}
D. Maclagan, B. Sturmfels {\it Introduction to Tropical Geometry} AMS, 2015

\bibitem[M06]{mikhicm} G. Mikhalkin, {\it Tropical Geometry and its applications}  Proceedings of the ICM 2006. Volume II: Invited lectures (2006) 827–852, arXiv:math/0601041





\bibitem[OS80]{os}
P. Orlik, L. Solomon, {\it Combinatories and Topology of Complements
of Hyperplanes}, lnvent. math. 56 (1980) 167-189

\bibitem[R98]{rybn}
G. Rybnikov, {\it On the fundamental group of the complement of a complex hyperplane arrangement}, Func.An.and Appl. 45 (2011) 137–148, arXiv:math/9805056

\bibitem[S93]{bs}
B. Shapiro, {\it The Mixed Hodge Structure of the Complement to an Arbitrary Arrangement of Affine Complex Hyperplanes is Pure}, Proc. AMS. 117 (1993) 931-933






\bibitem[T04]{tev}
J. Tevelev, {\it Compactifications of Subvarieties of Tori}, Amer. J. Math. 129 (2007) 1087-1104, arXiv:math/0412329

\bibitem[V76]{varch}
A. Varchenko, {\it Zeta-Function of Monodromy and Newton's Diagram}, Invent. math. 37 (1976) 253-262
\end{thebibliography}
\end{document}